\documentclass{birkjour}

\newtheorem{thm}{Theorem}[section]

\newtheorem{prop}[thm]{Proposition}
\theoremstyle{definition}
\newtheorem{defn}[thm]{Definition}
\newtheorem{assump}[thm]{Assumption}
\newtheorem{notation}[thm]{Notation}
\theoremstyle{remark}

\numberwithin{equation}{section}

\usepackage{enumerate}
\usepackage{url}

\begin{document}

\submitted{February 2, 2020}
\revised{April 16, 2020}
\accepted{April 17, 2020}


\title[Sufficient optimality conditions for delayed optimal control problems]%
{A survey on sufficient optimality conditions for delayed optimal control problems\footnote{
This is a preprint of the following paper: A.~P.~Lemos-Pai\~{a}o, C.~J.~Silva and D.~F.~M.~Torres, 
A survey on sufficient optimality conditions for delayed optimal control problems, published in 
Mathematical Modelling and Analysis of Infectious Diseases, edited by K.~Hattaf and H.~Dutta, 
Springer Nature Switzerland AG. Submitted 2/Feb/2020; revised 16/Apr/2020; accepted 17/Apr/2020.}}


\author[Lemos-Pai\~{a}o]{Ana~P.~Lemos-Pai\~{a}o}
\address{Center for Research and Development in Mathematics and Applications (CIDMA), 
Department of Mathematics, University of Aveiro, 3810-193 Aveiro, Portugal}
\email{anapaiao@ua.pt}


\author[Silva]{Cristiana~J.~Silva}
\address{Center for Research and Development in Mathematics and Applications (CIDMA), 
Department of Mathematics, University of Aveiro, 3810-193 Aveiro, Portugal}
\email{cjoaosilva@ua.pt}


\author[Torres]{Delfim~F.~M.~Torres}
\address{Center for Research and Development in Mathematics and Applications (CIDMA), 
Department of Mathematics, University of Aveiro, 3810-193 Aveiro, Portugal}
\email{delfim@ua.pt}


\begin{abstract}
The aim of this work is to make a survey on recent sufficient optimality conditions 
for optimal control problems with time delays in both state and control variables. 
The results are obtained by transforming delayed optimal control problems 
into equivalent non-delayed problems. Such approach allows to use standard 
theorems that ensure sufficient optimality conditions for non-delayed 
optimal control problems. Examples are given with the purpose 
to illustrate the results.
\end{abstract}


\thanks{This work is part of first author's Ph.D. project, 
carried out at the University of Aveiro under the support 
of the Portuguese Foundation for Science and Technology (FCT), 
fellowship PD/BD/114184/2016. It was also supported by FCT 
within project UIDB/04106/2020 (CIDMA)}



\maketitle


\section{Introduction}

The study of delayed systems, which can be optimized and controlled 
by a certain control function, has a long history and has been 
developed by many researchers (see, e.g., 
\cite{Banks,Bashier,Boccia1,Cacace,Elaiw,Friedman,Gollmann,Halanay,Oguztoreli,Stumpf} 
and references cited therein). Such systems can be called retarded, time-lag, 
or hereditary processes/optimal control problems. There are many applications 
of such systems, in diverse fields as Biology, 
Chemistry, Mechanics, Economy and Engineering 
(see, e.g., \cite{Bashier,Elaiw,Gollmann2,Ivanov,Klamka,Santos,Stumpf,Xu2,Xu}).
Dynamical systems with time delays, in both state and control variables, play
an important role in the modelling of real-life phenomena, in various fields
of applications (see \cite{Gollmann,Gollmann2}). For instance, in \cite{Rocha}
the incubation and pharmacological delays are modelled through the introduction
of time delays in both state and control variables. In \cite{Silva:2017},
Silva et al. introduce time delays in the state and control variables
for tuberculosis modelling. They represent the time delay on the diagnosis
and commencement of treatment of individuals with active tuberculosis infection
and the delays on the treatment of persistent latent individuals, due to clinical
and patient reasons. 

Delayed linear differential systems have also been investigated, their importance
being recognized both from a theoretical and practical points of view. For instance,
in \cite{Friedman}, Friedman considers linear hereditary processes and apply to them
Pontryagin's method, deriving necessary optimality conditions as well as existence and
uniqueness results. Analogously, in \cite{Oguztoreli}, delayed linear 
differential equations and optimal control problems involving this kind of systems
are studied. Since these first works, many researchers have devoted their attention
to linear quadratic optimal control problems with time delays 
(see, e.g., \cite{Cacace,Delfour,Eller,Khellat,Palanisamy}).
It turns out that for delayed linear quadratic optimal control problems
it is possible to provide an explicit formula for the optimal controls
(see \cite{Cacace,Khellat,Palanisamy}).

Delayed optimal control problems with differential systems, which are
linear both in state and control variables, have been studied in
\cite{Cacace,Chyung,Delfour,Eller,Khellat,Koepcke,Koivo,Lee1,Oguztoreli_2,Palanisamy}.
In \cite{Delfour,Koepcke,Palanisamy}, the system is delayed with respect to state 
and control variables. In \cite{Chyung,Oguztoreli_2}, the system only considers delays 
in the state variable. Chyung and Lee derive necessary and sufficient optimality 
conditions in \cite{Chyung}, while O\u{g}uzt\"{o}reli only proves necessary conditions 
in \cite{Oguztoreli_2}. Certain necessary conditions analysed by Chyung and Lee 
in \cite{Chyung} have been already derived in \cite{Kharatishvili_2,Pontryagin_et_all_1962,Popov}. 
However, the system considered in \cite{Chyung} is different from the previously studied
hereditary systems, which do not require a initial function of state.
In \cite{Eller}, Eller et al. derive a sufficient condition for a control
to be optimal for certain problems with time delay. The problems studied
by Eller et al. and Khellat in \cite{Eller} and \cite{Khellat}, respectively, 
consider only one constant lag in the state. The research done by Lee in \cite{Lee1}
is different from that of the current work (more specifically from that 
of Section~\ref{section-delayed-linear}), because in \cite{Lee1} the aim is 
to minimize a cost functional, which does not consider delays, subject to a linear 
differential system (with respect to state and control variables) 
and to another constraint. In their differential system, the state variable 
depends on a constant and a fixed delay, while the control variable depends 
on a constant lag, which is not specified a priori. Note that the
differential system of the problem considered in \cite{Koivo} is similar
to the one of \cite{Lee1}. Although Banks has studied delayed non-linear 
problems without lags in the control, he has also analysed problems that
are linear and delayed with respect to control (see \cite{Banks}). 
Later, in 2010, Carlier and Tahraoui investigated optimal control problems with 
a unique delay in the state (see \cite{Carlier}). 
In 2012 and 2013, Frederico and Torres devoted their attention to optimal 
control problems that only contain delays in the state variables 
and the dependence on the control is linear (see \cite{Frederico3,Frederico1}).  
Recently, Cacace et al. studied optimal control problems that involve linear 
differential systems with variable delays only in the control (see \cite{Cacace}).

The problems analysed in the current work are different from those considered
in the mentioned papers. In Section~\ref{section-delayed-linear}, the optimal 
control problems involve differential systems that are linear with respect to state, 
but not with respect to the control. In Section~\ref{section-delayed-non-linear}, 
we study optimal control problems with non-linear differential systems.
Furthermore, in both Sections~\ref{section-delayed-linear} and 
\ref{section-delayed-non-linear}, we consider a constant time delay 
in the state and another one in the control.
These two delays are, in general, not equal.

In \cite{Hughes}, Hughes firstly consider variational problems with only
one constant lag and derive various necessary and a sufficient optimality
conditions for them. The variational problems in \cite{Hughes} can easily
be transformed into control problems with only one constant delay (see, e.g.,
\cite[p.~53--54]{Paiao}). Hughes also investigates an optimality condition for
a control problem with a constant delay, which is the same for state and control.
The problems analysed by Chan and Yung in \cite{Chan}, and by Sabbagh in \cite{Sabbagh}, 
are similar to the first problems studied by Hughes in \cite{Hughes}. 
Therefore, the problems investigated in \cite{Chan,Hughes,Sabbagh} are different 
from the problems studied by us here, because in the present work the state delay 
is not necessarily equal to the control delay.
The problems considered in \cite{Hughes,Sabbagh} are also considered in \cite{Palm}
by Palm and Schmitendorf. For such problems, they derive two conjugate-point conditions,
which are not equivalent. Note that their conditions are only necessary
and do not give a set of sufficient conditions (see \cite{Palm}).
Recent results include Noether type theorems for problems of the calculus 
of variations with time delays (see \cite{Frederico2,Malinowska,Santos2}),
necessary optimality conditions for quantum (see \cite{Frederico1})
and Herglotz variational problems with time delays (see \cite{Santos3,Santos}),
as well as delayed optimal control problems with integer 
(see \cite{Benharrat,Bokov,Frederico3}) and non-integer 
(fractional order) dynamics (see \cite{Debbouche1,Debbouche2}).
Applications of such theoretical results are found in Biology 
and other Natural Sciences, e.g., in tuberculosis (see \cite{Silva:2017})
and HIV (see \cite{Rocha,Rodrigues}). 

In \cite{Jacobs}, Jacobs and Kao investigate delayed problems
that consist to minimize a cost functional without delays subject
to a differential system defined by a non-linear function
with a delay in state and another one in the control. Similarly to
our problems, these delays do not have to be equal. In contrast,
all types of cost functionals considered in our work also contain time delays.
Therefore, we study here problems that are more general than the one considered 
in \cite{Jacobs}. Jacobs and Kao transform the problem using a Lagrange-multiplier 
technique and prove a regularity result in the form of a controllability condition,
as well as some necessary optimality conditions. Then, in some special
restricted cases, they prove existence, uniqueness, and sufficient conditions.
Such restricted problems consider a differential system that is linear
in state and in control variables. Thus, the sufficient conditions
of \cite{Jacobs} are derived for problems that are less general than ours.

As it is well-known, and as Hwang and Bien write in \cite{Hwang}, 
many researchers have directed their efforts to seek sufficient 
optimality conditions for control problems with delays 
(see, e.g., \cite{Chyung,Eller,Hughes,Jacobs,Lee3,Schmitendorf}). 
Therefore, it is not a surprise that there are authors that already 
proved some sufficient optimality conditions for delayed optimal 
control problems similar but, nevertheless, different from ours. 
In what respects to research done in \cite{Chyung,Eller,Hughes,Jacobs},
we have already seen why they are different. The delayed optimal control 
problems analysed by Schmitendorf in \cite{Schmitendorf} have a cost 
functional and a differential system that are more general than ours. 
However, in \cite{Schmitendorf} the control takes its values 
in all $\mathbb{R}^m$, while in the present work
the control values belong to a set $U\subseteq\mathbb{R}^m$,
$m\in\mathbb{N}$. In \cite{Lee3}, Lee and Yung study a problem
that is similar to the one considered in \cite{Schmitendorf},
where the control belongs to a subset of $\mathbb{R}^m$, as we consider here.
First and second-order sufficient conditions are shown in \cite{Lee3}.
Nevertheless, the conditions of \cite{Lee3} are not constructive
and not practical for the computation of the optimal solution.
Indeed, as hypothesis, it is assumed the existence of a symmetric matrix
under some conditions, for which is not given a method to calculate
its expression. Another similar problem to ours is studied by Bokov
in \cite{Bokov}, in order to arise a necessary optimality condition
in an explicit form. Moreover, a solution to the problem with
infinite time horizon is given in \cite{Bokov}. In contrast, in the 
present work we are interested to derive sufficient optimality conditions.
In \cite{Hwang}, Hwang and Bien prove a sufficient condition for problems
involving a differential affine time delay system with the same time delay 
for the state and the control. The differential systems considered in 
the present work are more general. In 1996, Lee and Yung, considering 
functions that do not have to be convex, derived various first and 
second-order sufficient conditions for non-linear optimal control problems 
with only a constant delay in the state (see \cite{Lee2}). 
Their class of problems is obviously different from our.
In particular, we consider delays for both state 
and control variables. As in \cite{Chan,Lee3}, second-order sufficient
conditions are shown to be related to the existence of solutions
of a Riccati-type matrix differential inequality.

Optimal control problems with multiple delays have also been investigated.
In \cite{Halanay}, Halanay derive necessary conditions for some optimal
control problems with various time lags in state and control variables,
using the abstract multiplier rule of Hestenes (see \cite{Hestenes}).
In \cite{Halanay}, all delays related to state are equal to each other
and the same happens with the delays associated to the control.
Note that the results of \cite{Friedman,Kharatishvili_2} are obtained
as particular cases of problems considered in \cite{Halanay}. Later, 
in 1973, a necessary condition is derived for an optimal control problem
that involves multiple constant lags only in the control. This delayed
dependence occurs both in the cost functional and in the differential
system, which is defined by a non-linear function (see \cite{Soliman}).
In \cite{Kharatishvili_3}, Harati\v{s}vili and Tadumadze prove the existence
of an optimal solution and a necessary condition for optimal control systems
with multiple variable time lags in the state and multiple variable
commensurable time delays in the control. Later, an optimal control
problem where the state variable is solution of an integral equation
with multiple delays, both on state and control variables,
is studied by Bakke in \cite{Bakke}. Furthermore, necessary conditions 
and Hamilton--Jacobi equations are derived. In 2006, Basin 
and Rodriguez-Gonzalez proved a necessary and a sufficient optimality 
condition for a problem that consists to minimize 
a quadratic cost functional subject to a linear system 
with multiple time delays in the control variable (see \cite{Basin}). 
In their work, they begin by deriving a necessary condition through 
Pontryagin's Maximum Principle. Afterwards, sufficiency is proved 
by verifying if the candidate found, through the Maximum Principle, 
satisfies the Hamilton--Jacobi--Bellman equation. Although Basin 
and Rodriguez-Gonzalez consider multiple time delays, the dependence 
of the state and control in the differential system is linear. 
In our current work, the dependence of the control, 
in the differential systems, is, in general, non-linear.
In 2013, Boccia et al. derived necessary conditions
for a free end-time optimal control problem subject to a non-linear
differential system with multiple delays in the state (see \cite{Boccia1}).
The control variable is not influenced by time lags in \cite{Boccia1}.
Recently, in 2017, Boccia and Vinter obtained necessary conditions for
a fixed end-time problem with a constant and unique delay for all variables,
as well as free end-time problems without control delays (see \cite{Boccia2}).

As Guinn wrote in \cite{Guinn}, the classical methods of obtaining necessary 
conditions for retarded optimal control problems (used, for instance, 
by Halanay in \cite{Halanay}, Harati\v{s}vili in \cite{Kharatishvili}
and O\v{g}uzt\"{o}reli in \cite{Oguztoreli}) require complicated
and extensive proofs (see, e.g.,
\cite{Banks,Friedman,Halanay,Kharatishvili,Oguztoreli}). In $1976$,
Guinn proposed a method whereby we can reduce some specific time-lag
optimal control problems to equivalent and augmented optimal control
problems without delays (see \cite{Guinn}). By reducing delayed optimal control
problems into non-delayed ones, we can then use well-known theorems,
applicable for optimal control problems without delays, to derive
desired optimality conditions for delayed problems (see \cite{Guinn}).
In \cite{Guinn}, Guinn study specific optimal control problems with 
a constant delay in state and control variables. These two delays are equal.
Later, in $2009$, G\"{o}llmann et al. studied optimal control problems
with a constant delay in state and control variables subject to
mixed control-state inequality constraints (see \cite{Gollmann}). 
In that research, the delays do not have to be equal. For technical reasons,
the authors need to assume that the ratio between these two time delays
is a rational number (see \cite{Gollmann}). In \cite{Gollmann}, the method used
by Guinn in \cite{Guinn} is generalized and, consequently, a non-delayed optimal
control problem is again obtained. Pontryagin's Minimum Principle,
for non-delayed control problems with mixed state-control constraints,
is used and first-order necessary optimality conditions are derived
for retarded problems (see \cite{Gollmann}). Furthermore, G\"{o}llmann et al. 
discuss the Euler discretization of the retarded problem and some analytical 
examples versus correspondent numerical solutions are given. 
For more on numerical methods, for solving applied optimal control problems 
of systems governed by delay differential equations, see \cite{MR3318705,ref2,MR3723487}. 
Later, in 2014, G\"{o}llmann and Maurer generalized the research mentioned before, 
by studying optimal control problems with multiple and constant time delays 
in state and control, involving mixed state-control inequality constraints 
(see \cite{Gollmann2}). Again, necessary optimality conditions are derived 
(see \cite{Gollmann2}). Note that the works \cite{Gollmann,Gollmann2,Guinn,Halanay}
consider delayed non-linear differential systems.

In Section~\ref{section-delayed-linear}, we consider optimal control problems
that consist to minimize a delayed non-linear cost functional
subject to a delayed differential system that is linear
with respect to state, but not with respect to control.
Note that the cost functional does not have to be quadratic,
but it satisfies some continuity and convexity assumptions.
In Section~\ref{section-delayed-non-linear}, we consider optimal control problems
that consist to minimize a delayed non-linear cost functional
subject to a delayed non-linear differential system.
In both Sections~\ref{section-delayed-linear} and \ref{section-delayed-non-linear},
the delay in the state is the same for the cost functional
and for the differential system. The same happens with
the time lag of the control variable.
Analogously to G\"{o}llmann et al. in \cite{Gollmann}, we ensure 
the Commensurability Assumption between the, 
possibly different, delays of state and control variables. 
The proofs of our sufficient optimality conditions consider the technique 
proposed by Guinn in \cite{Guinn} and used by G\"{o}llmann et al. 
in \cite{Gollmann,Gollmann2} (see \cite{Lemos-Paiao1,Lemos-Paiao2}). 
As we have already mentioned before, the technique consists 
to transform a delayed optimal control problem into 
an equivalent non-delayed optimal control problem. 
After doing such transformation, one can apply well-known results 
for non-delayed optimal control problems and then return 
to the initial delayed problem. Here we restrict ourselves to delayed problems 
with deterministic controls. For the stochastic case, we refer the reader 
to \cite{Frederico2,Fuhrman,Goldys,Ivanov,Larssen}.

This work is organised as follows. We begin by recalling the 
Commensurability Assumption, introduced by G\"{o}llmann et al. in \cite{Gollmann}, 
and by defining some needed notations, in Section~\ref{sec:assumptions}.
In Section~\ref{subsec:def:linear}, we define a state-linear optimal control
problem with constant time delays in state and control variables.
Then, in Section~\ref{subsec:theo:linear}, we present a sufficient optimality
condition associated with the problem stated in Section~\ref{subsec:def:linear}.
A concrete example is solved in detail in Section~\ref{subsec:example:linear},
with the purpose to illustrate Theorem~\ref{Cap4:theo_OCP-LD} of 
Section~\ref{subsec:theo:linear}. In Section~\ref{subsec:theo:non-linear}, 
we present a sufficient optimality condition associated with the non-linear 
optimal control problem with time lags both in state and control variables, 
defined in Section~\ref{subsec:def:non-linear}. An example that 
illustrates the obtained theoretical result -- Theorem~\ref{Cap4:theo_suf_NLD} 
of Section~\ref{subsec:theo:non-linear} -- is given.  We end with some 
conclusions, in Section~\ref{conclusion}.


\section{Commensurability assumption and notations}
\label{sec:assumptions}

In this section, we recall the Commensurability Assumption 
introduced by G\"{o}llmann et al. in \cite{Gollmann}.

\begin{assump}[See Assumption~4.1 of \cite{Gollmann}]
\label{assump:commensurability}
We consider $r,s\geq 0$, not simultaneously equal to zero, 
and commensurable, that is, 
$$
(r,s)\neq(0,0)
$$ 
and
\begin{equation*}
\frac{r}{s}\in\mathbb{Q}\ \text{ for }\ s>0
\ \text{ or }\ \frac{s}{r}\in\mathbb{Q}\ \text{ for }\ r>0.
\end{equation*}
\end{assump}

Actually, Commensurability Assumption~\ref{assump:commensurability} 
holds for any couple of rational numbers $(r, s)$ for which 
at least one number is non-zero (see \cite{Gollmann}).

With the purpose to simplify the writing, we introduce some notations. 

\begin{notation}
\label{notation:time-delays}
We define $t_{\tau}$, $t^{\tau}$ 
and $t_{\tau_1}^{\tau_2}$ as follows:
$$
t_{\tau}=t-\tau,\quad t^{\tau}=t+\tau 
\quad \text{and} \quad t^{\tau_2}_{\tau_1}=t-\tau_1+\tau_2
$$	
for time delays $\tau,\tau_1,\tau_2\in\{r,s\}$
and for all $t\in[a,b]$.
\end{notation}

\begin{notation}
\label{notation_NLD}
Let 
$x_a=x(a)=\varphi(a)$
and
$x_r(t)=\big (x(t),x(t-r)\big )$.
Moreover, we define the operators
$[\cdot,\cdot]_r$ and $\langle\cdot,\cdot\rangle_r$ by
$[x,\zeta]_r(t) := \Big (t,x_r(t),\zeta\big (t,x_r(t)\big )\Big )$
and
$\langle x,\zeta\rangle_r(t) := \Big (t,x_r(t),\zeta\big (t,x(t)\big )\Big )$,
respectively.
\end{notation}

While Notation~\ref{notation:time-delays} is used 
in Sections~\ref{section-delayed-linear} and \ref{section-delayed-non-linear}, 
Notation~\ref{notation_NLD} is only considered 
in Section~\ref{section-delayed-non-linear}.


\section{Delayed state-linear optimal control problem}
\label{section-delayed-linear}

This section is devoted to state-linear optimal control problems
with constant time delays in state and control variables.
We make a survey on a sufficient optimality condition for this type of problems.
Its proof, and more details associated with the contents of the current section, 
can be found in \cite{Lemos-Paiao1}. To finish this section, 
an illustrative example is given.  


\subsection{Statement of the optimal control problem}
\label{subsec:def:linear}

We start by defining a delayed state-linear optimal control problem.

\begin{defn}
\label{Cap4:def:OCP_LD}
Consider that $r \geq 0$ and $s \geq 0$ are constant
time delays associated with the state and control variables, respectively.
We assume that $(r, s) \neq (0, 0)$.
A non-autonomous state-linear optimal control problem
(OCP$_\text{LD}$) with time delays and
with a fixed initial state, on a fixed finite time interval $[a,b]$, consists in
\begin{equation*}
\min\ \ C_{LD}\big (x(\cdot),u(\cdot)\big )=\int_{a}^{b}f^0_x\big (t,x(t),x(t-r)\big )+f^0_u\big (t,u(t),u(t-s)\big )dt
\end{equation*}
subject to the delayed differential system
\begin{equation}
\label{eq:Delay:linear:CS}
\dot{x}(t)=A(t)x(t)+A_D(t)x(t-r)+g\big (t,u(t)\big )+g_D\big (t,u(t-s)\big )
\end{equation}
with the following initial conditions:
\begin{equation}
\label{eq:InitCond:delayControl}
\begin{split}
x(t)&=\varphi(t),\ t\in[a-r,a],\\
u(t)&=\psi(t),\ t\in[a-s,a[,
\end{split}
\end{equation}
where
\begin{enumerate}[i.]
\item the state trajectory is 
$x(t)\in\mathbb{R}^n$ 
for each $t\in[a-r,b]$;
\item the control is 
$u(t)\in U\subseteq\mathbb{R}^m$ 
for each $t\in[a-s,b]$;
\item $A(t)$ and $A_D(t)$ are real $n \times n$ matrices 
for each $t\in[a,b]$.
\end{enumerate}
\end{defn}

Next we define admissible pair for (OCP$_\text{LD}$).
	
\begin{defn}
We say that $\big (x(\cdot),u(\cdot)\big )$ is an admissible pair  
for (OCP$_\text{LD}$) if it respects the following conditions:
\begin{enumerate}[i.]
\item $\big (x(\cdot),u(\cdot)\big )\in W^{1,\infty}([a-r,b],
\mathbb{R}^n)\times L^{\infty}([a-s,b],\mathbb{R}^m)$,
where $W^{1,\infty}$ is the space of Lipschitz functions;
\item $\big (x(\cdot),u(\cdot)\big )$ satisfies the conditions \eqref{eq:Delay:linear:CS} 
and \eqref{eq:InitCond:delayControl};
\item $\big (x(t),u(t)\big )\in \mathbb{R}^n\times U$ for all $t\in[a,b]$.
\end{enumerate}
\end{defn}
	

\subsection{Main result}
\label{subsec:theo:linear}

In what follows, we consider that the time delays 
$r$ and $s$ respect the Commensurability Assumption~\ref{assump:commensurability} 
and we use Notation~\ref{notation:time-delays}.
	
The following theorem supplies a sufficient optimality 
condition associated with (OCP$_\text{LD}$) 
(see Definition~\ref{Cap4:def:OCP_LD}).
Such result generalizes Theorem~5 in Chapter 5.2 of \cite{Lee}.
	
\begin{thm}
\label{Cap4:theo_OCP-LD}
Consider \emph{(OCP$_\text{LD}$)} and assume that
\begin{enumerate}[i.]
			
\item functions $f^0_x$, $\partial_2f^0_x$, $\partial_3f^0_x$, $f^0_u$,
$g$, $g_D$, $A$ and $A_D$ are continuous for all their arguments;
			
\item $f^0_x(t,x,x_r)$ is a convex function in $(x,x_r)\in\mathbb{R}^{2n}$
for each $t\in[a,b]$;
		
\item for almost all $t\in[a,b]$, $u^*$ is a control with response $x^*$
that satisfies the \emph{maximality condition}
\begin{equation*}
\begin{split}
&\ \max_{u\in U}\big \{H_D^1\big (t,x^*(t),x^*(t_r),u,u^*(t_s),\eta(t)\big )\\
&\ + H_D^0\big (t^s,x^*(t^s),x^*(t^s_r),u^*(t^s),u,\eta(t^s)\big )\chi_{[a,b-s]}(t)\big \}\\
= &\ H_D^1\big (t,x^*(t),x^*(t_r),u^*(t),u^*(t_s),\eta(t)\big )\\
&\ + H_D^0\big (t^s,x^*(t^s),x^*(t^s_r),u^*(t^s),u^*(t),\eta(t^s)\big )\chi_{[a,b-s]}(t),
\end{split}
\end{equation*}
where
\begin{equation*}
\begin{split}
H_D^p(t,x,y,u,v,\eta)=&-\Big [f^0_x(t,x,y)+f^0_u(t,u,v)\Big ]\\
&+\eta \Big [A(t)x+A_D(t)y+pg(t,u)+(1-p)g_D(t,v)\Big ]
\end{split}
\end{equation*}
for $p\in\{0,1\}$ and $\eta(t)$ is any non-trivial solution
of the adjoint system
\begin{equation*}
\begin{split}
\dot{\eta}(t)=&\ \partial_2f^0_x\big (t,x^*(t),x^*(t_r)\big )+\partial_3
f^0_x\big (t^r,x^*(t^r),x^*(t)\big )\chi_{[a,b-r]}(t)\\
&-\eta(t)A(t)-\eta(t^r)A_D(t^r)\chi_{[a,b-r]}(t)
\end{split}
\end{equation*}
that satisfies the transversality condition
$\eta(b)=
[
\begin{matrix}
0 & \cdots & 0
\end{matrix}
]_{1\times n}
$.
\end{enumerate}
Then, $\big (x^*(\cdot),u^*(\cdot)\big )$ is an optimal solution 
of \emph{(OCP$_\text{LD}$)} that leads to the minimal cost 
$C_{LD}\big (x^*(\cdot),u^*(\cdot)\big )$.
\end{thm}

The detailed proof of Theorem~\ref{Cap4:theo_OCP-LD} 
can be found in \cite{Lemos-Paiao1}.


\subsection{An illustrative example}
\label{subsec:example:linear}

In this section we provide an illustrative example 
associated with Theorem~\ref{Cap4:theo_OCP-LD}.
	
\begin{figure}[htp]
\vspace*{-10pt}\begin{center}
\includegraphics[scale=0.37]{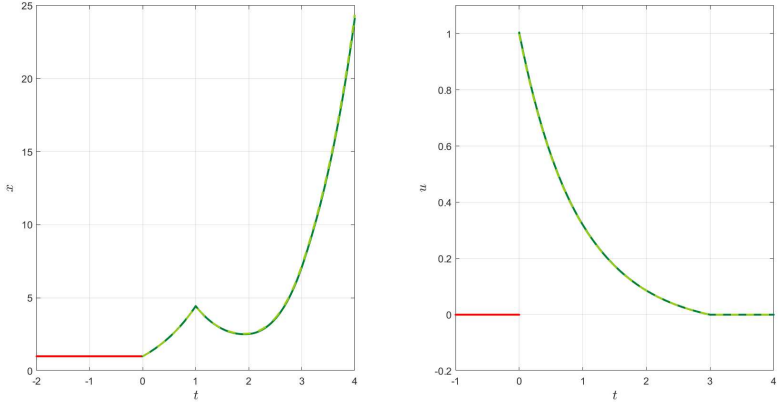}
\caption{Optimal solution of problem \eqref{example}: 
red line -- initial data;
dark green line -- analytical solution;
dashed light green line -- numerical solution.}
\label{Cap4:plot-exemplo}
\end{center}
\end{figure}

Let us consider the delayed state-linear optimal control problem given by
\begin{equation}
\label{example}
\begin{split}
\min\ \ \ & C_{LD}\big (x(\cdot),u(\cdot)\big )=\int_{0}^{4}x(t)+100u^2(t)dt\\
\text{s.t.}\ \ & \dot{x}(t)=x(t)+x(t-2)-10u(t-1),\\
& x(t)=1,\ t\in[-2,0],\\
& u(t)=0,\ t\in[-1,0[,
\end{split}
\end{equation}
where $u(t)\in U=\mathbb{R}$ for each $t\in[-1,4]$. Thus,
we have that $n=m=1$, $a=0$, $b=4$, $r=2$, $s=1$, 
$f^0_x\big (t,x(t),x(t-2)\big )=x(t)$,
$f^0_u\big (t,u(t),u(t-1)\big )=100u^2(t)$, $A(t)=A_D(t)=1$, $g\big (t,u(t)\big )=0$
and $g_D\big (t,u(t-1)\big )=-10u(t-1)$. Note that our functions respect
hypothesis $i$ and $ii$ of Theorem~\ref{Cap4:theo_OCP-LD}. Let $\bar{u}$
be an admissible control of problem~\eqref{example} and let us maximize function
\begin{equation*}
\begin{split}
&-f^0_u\big (t,u,\bar{u}(t-1)\big )+\eta(t)g(t,u)\\
&+\big [-f^0_u\big(t+1,\bar{u}(t+1),u\big)+\eta(t+1)g_D(t+1,u)\big ]\chi_{[0,3]}(t)\\
=&-100u^2+\big [-100\bar{u}^2(t+1)-10\eta(t+1)u\big ]\chi_{[0,3]}(t)\\
=&
\begin{cases}
-100u^2-10\eta(t+1)u-100\bar{u}^2(t+1),\ & t\in[0,3]\\
-100u^2,\ & t\in\ ]3,4]
\end{cases}
\end{split}
\end{equation*}
with respect to $u\in\mathbb{R}$. We obtain
\begin{equation*}
u(t)=-\frac{\eta(t+1)}{20}
\end{equation*}
for $t\in[0,3]$ and $u(t)=0$ for $t\in\ ]3,4]$.
Furthermore, we know that $\eta(t)$ is any non-trivial solution of
\begin{equation*}
\begin{split}
\dot{\eta}(t)=&\ \partial_2f^0\big(t,x(t),x(t-2)\big)
+\partial_3f^0\big(t+2,x(t+2),x(t)\big)\chi_{[0,2]}(t)\\
&-\eta(t)A(t)-\eta(t+2)A_D(t+2)\chi_{[0,2]}(t)\\
\Leftrightarrow
\dot{\eta}(t)=&\ 1-\eta(t)-\eta(t+2)\chi_{[0,2]}(t)
=
\begin{cases}
1-\eta(t)-\eta(t+2),\ & t\in[0,2]\\
1-\eta(t),\ & t\in\ ]2,4]
\end{cases}
\end{split}
\end{equation*}
that satisfies the transversality condition $\eta(4)=0$.
The adjoint system is given by
\begin{equation}
\begin{split}
\label{adjoint_system}
\begin{cases}
\dot{\eta}(t)=
\begin{cases}
1-\eta(t)-\eta(t+2),\ & t\in[0,2]\\
1-\eta(t),\ & t\in\ ]2,4]\\
\end{cases}\\
\eta(4)=0.
\end{cases}
\end{split}
\end{equation}
For $t\in\ ]2,4]$, the solution of differential equation
\begin{equation*}
\begin{cases}
\dot{\eta}(t)
= 1-\eta(t)\\
\eta(4)=0
\end{cases}
\end{equation*}
is given by $\eta(t)=1-e^{4-t}$.
Knowing $\eta(t)$, $t\in\ ]2,4]$, and attending to the continuity
of function $\eta$ for all $t\in[0,4]$, we can determine $\eta(t)$
for $t\in[0,2]$ solving the differential equation
\begin{equation*}
\begin{cases}
\dot{\eta}(t)
=1-\eta(t)-\eta(t+2)\\
\eta(2)=1-e^{4-2}=1-e^2
\end{cases}
\end{equation*}
for $t\in[0,2]$. Therefore, we have that $\eta(t)=e^{2-t}\left(t-e^2-1\right)$ for $t\in[0,2]$.
Consequently, the solution of the adjoint system~\eqref{adjoint_system} is given by
\begin{equation*}
\begin{split}
\eta(t)
=
\begin{cases}
e^{2-t}\left(t-e^2-1\right),\ & t\in [0,2]\\
1-e^{4-t},\ & t\in\ ]2,4].
\end{cases}
\end{split}
\end{equation*}
So, the control is given by
\begin{equation}
\begin{split}
\label{optimalcontrol}
u(t)=\frac{1}{20}
\begin{cases}
0,\ & t\in[-1,0[\\
e^{3-t}-e^{1-t}t,\ & t\in[0,1[\\
e^{3-t}-1,\ & t\in[1,3]\\
0,\ & t\in\ ]3,4].
\end{cases}
\end{split}
\end{equation}
Knowing the control, we can determine
the state by solving the differential equation
\begin{equation*}
\begin{cases}
\dot{x}(t)=x(t)+x(t-2)-10u(t-1)\\
x(t)=1,\ t\in[-2,0].
\end{cases}
\end{equation*}
The state solution is given by
\begin{equation}
\label{optimalstate}
\small
x(t) = 
\begin{cases}
1,& t\in[-2,0]\\[0.8em]
-1+2e^t,& t\in\ ]0,1]\\[0.8em]
\displaystyle\frac{\left(e^2+2e^4-2e^2t\right)e^{-t}-8+\left(17-2e^2\right)e^t}{8},& t\in\ ]1,2]\\[0.8em]
\displaystyle\frac{2e^{4-t}+4+\left(-47e^{-2}+17-2e^2+16e^{-2}t\right)e^t}{8},& t\in\ ]2,3]\\[0.8em]
\displaystyle\frac{\left(-e^6+e^4t\right)e^{-t}+4+\left(-51e^{-2}+24-2e^2+17e^{-2}t-2t\right)e^t}{8},& t\in\ ]3,4].
\end{cases}
\end{equation}
Such analytical expressions can be obtained with the help of a modern computer
algebra system. We have used \textsf{Mathematica}.
In Figure~\ref{Cap4:plot-exemplo}, we observe that the numerical
solutions for control and state, obtained using \textsf{AMPL} \cite{AMPL} and
\textsf{IPOPT} \cite{IPOPT}, are in agreement with their analytical solutions, given by
\eqref{optimalcontrol} and \eqref{optimalstate}, respectively. The numerical solutions
were obtained using Euler's forward difference method in \textsf{AMPL} and \textsf{IPOPT},
dividing the interval of time $[0,4]$ into 2000 subintervals. The minimal cost is
$$
\frac{23+e^2+34e^4-2e^6}{16}\simeq 67.491786.
$$


\section{Delayed non-linear optimal control problem}
\label{section-delayed-non-linear}

This section is devoted to non-linear optimal control problems 
with constant time delays in state and control variables.
We make a survey on a sufficient optimality condition for this type of problems.
Its proof, and more details associated with the contents of the current section, 
can be found in \cite{Lemos-Paiao2}. We finish this section with 
an illustrative example.


\subsection{Statement of the optimal control problem}
\label{subsec:def:non-linear}
	
We start by defining the delayed 
non-linear optimal control problem.

\begin{defn}
\label{Cap3:def:OCP_D}
Consider that $r \geq 0$ and $s \geq 0$ are constant
time delays associated with the state and control variables, respectively.
A non-autonomous optimal control problem with constant time delays and
with a fixed initial state, on a fixed finite time interval $[a,b]$, 
is denoted by (OCP$_\text{D}$) and consists in
\begin{equation*}
\min\ \ C_D\big (x(\cdot), u(\cdot)\big )=g^0\big (x(b)\big ) + \int_{a}^{b}f^0\big (t,x(t),x(t-r),u(t),u(t-s)\big )dt
\end{equation*}
subject to the delayed differential system
\begin{equation}
\label{eq_dif_syst_D}
\dot{x}(t)=f\big (t,x(t),x(t-r),u(t),u(t-s)\big )
\quad \text{for a.e.} \quad t\in[a,b]
\end{equation}
with initial and final conditions
\begin{equation}
\begin{split}
\label{eq_init_cond_D}
x(t) &= \varphi(t),\quad t\in [a-r-s,a]\subset\mathbb{R},\\
u(t) &= \psi(t),\quad t\in[a-s,a[,\\
x(b) &\in \Pi\subseteq\mathbb{R}^n;
\end{split}
\end{equation}
where
\begin{enumerate}[i.]
\item the state trajectory is 
$x(t)\in\mathbb{R}^n$ 
for all $t\in [a-r-s,b]$;
\item the control is 
$u(t)\in U\subseteq\mathbb{R}^m$ 
for all $t\in[a-s,b]$;
\item $f=
\left[\begin{matrix}
f_1 & \cdots & f_n
\end{matrix}\right]^T$.
\end{enumerate}		
\end{defn}

Next we define admissible pair for (OCP$_\text{D}$).

\begin{defn}
We say that $\big (x(\cdot),u(\cdot)\big )$ is an admissible pair 
for (OCP$_\text{D}$) if it
respects the following conditions:
\begin{enumerate}[i.]
\item $\big (x(\cdot), u(\cdot)\big )\in W^{1,\infty}\big ([a-r-s,b],
\mathbb{R}^n\big )\times L^{\infty}\big ([a-s,b],\mathbb{R}^m\big )$;
\item $\big (x(\cdot), u(\cdot)\big )$ satisfies 
conditions~\eqref{eq_dif_syst_D} and \eqref{eq_init_cond_D};
\item $\big (x(t),u(t)\big )\in \mathbb{R}^n\times U$ for all $t\in[a,b]$.
\end{enumerate}
\end{defn}


\subsection{Main result}
\label{subsec:theo:non-linear}

In what follows we also consider that the time delays $r$ and $s$ respect 
Commensurability Assumption~\ref{assump:commensurability}.
Moreover, here we use Notations~\ref{notation:time-delays} and \ref{notation_NLD}.

The following theorem provides a sufficient optimality condition associated with  
(OCP$_\text{D}$) (see Definition~\ref{Cap3:def:OCP_D}).
Such result generalizes Theorem~7 in Chapter 5.2 of \cite{Lee}.

\begin{thm}
\label{Cap4:theo_suf_NLD}
Consider \emph{(OCP$_\text{D}$)}.
Let the interval $[a,b]$ be divided into $N \in\mathbb{N}$ subintervals of amplitude 
$h = \frac{b-a}{N}>0$ and suppose that the functions $g^0$, $f^0$ and $f$ are of 
class $\mathcal{C}^1$ with respect to all their arguments. 
Assume there exists a $\mathcal{C}^1\left(\mathbb{R}^{1+3n},\mathbb{R}^m\right)$ 
feedback control $u^*\Big (t,x_r(t),\eta\big (t,x_r(t)\big )\Big ) = u^*[x,\eta]_r(t)$ such that
\begin{equation*}
\begin{split}
&\ \max_{u\in U} \Big \{ H\Big (t,x_r(t),u,u^*[x,\eta]_r(t_s),\eta\big (t,x_r(t)\big )\Big )\\
&\ + H\Big (t^s,x_r(t^s),u^*[x,\eta]_r(t^s),
u,\eta\big (t^s,x_r(t^s)\big )\Big )\chi_{[a,b-s]}(t)\Big \}\\
= &\ H\Big (t,x_r(t),u^*[x,\eta]_r(t),u^*[x,\eta]_r(t_s),\eta\big (t,x_r(t)\big )\Big )\\
&\ + H\Big (t^s,x_r(t^s),u^*[x,\eta]_r(t^s),
u^*[x,\eta]_r(t),\eta\big (t^s,x_r(t^s)\big )\Big )\chi_{[a,b-s]}(t)\\
=: &\ H^0[x,\eta]_r(t)+H^0[x,\eta]_r(t^s)\chi_{[a,b-s]}(t)
\end{split}
\end{equation*}
for all $t\in[a,b]$, where
\begin{equation*}
H(t,x,y,u,v,\eta)=-f^0(t,x,y,u,v)+ \eta f(t,x,y,u,v).
\end{equation*}
Furthermore, let $I_i=[a+hi,a+h(i+1)]$, 
$i=0, \ldots, N-1$, and suppose that function 
$S\big (t,x(t)\big )\in\mathcal{C}^2\left(\mathbb{R}^{1+n},
\mathbb{R}\right)$, $t\in[a,b]$, 
is a solution of equation
\begin{equation}
\label{eq_HJ_main_theo}
\begin{split}
&\ \partial_1S\big (t,x(t)\big )+\sum_{i=0}^{N-1}\Big\{-f^0\big (t,
x_r(t),u^*\langle x,\partial_2 S\rangle_r(t),
u^*\langle x,\partial_2 S\rangle_r(t_s)\big )\\
&+\partial_2S\big (t,x(t)\big )f\big (t,x_r(t),u^*\langle x,\partial_2 S\rangle_r(t),
u^*\langle x,\partial_2 S\rangle_r(t_s)\big )\Big\}\chi_{I_i}(t)=0
\end{split}
\end{equation}
with $S\big (b,x(b)\big )=-g^0\big (x(b)\big )$, $x(b)\in \Pi$.
Finally, consider that the control law
$$
u^*\Big (t,x_r(t),\partial_2S\big (t,x(t)\big )\Big ) 
= u^*\langle x,\partial_2 S\rangle_r(t),
\quad t\in[a,b],
$$
determines a response $\tilde{x}(t)$ 
steering $(a,x_a)$ to $(b,\Pi)$. Then, 
$$
\tilde{u}(t)
=u^*\big (t,\tilde{x}(t),\tilde{x}(t-r),\partial_2S(t,\tilde{x}(t)\big )
$$
is an optimal control of \emph{(OCP$_\text{D}$)} that leads to the minimal cost 
$$
C_D\big (\tilde{x}(\cdot), \tilde{u}(\cdot)\big )=-S(a,x_a).
$$
\end{thm}

The detailed proof of Theorem~\ref{Cap4:theo_suf_NLD} 
can be found in \cite{Lemos-Paiao2}.


\subsection{An illustrative example}
\label{subsec:example:non-linear}

In this section we provide an example of application 
of Theorem~\ref{Cap4:theo_suf_NLD}.
	
Let us consider the following delayed non-linear optimal control problem
studied by G\"{o}llmann et al. in \cite{Gollmann}:
\begin{equation}
\label{example1}
\begin{split}
\min\ \ \ & C_D\big (x(\cdot),u(\cdot)\big )=\int_{0}^{3} x^2(t)+u^2(t)dt,\\
\text{s.t.}\ \ & \dot{x}(t)=x(t-1)\ u(t-2),\\
& x(t)=1,\quad t\in[-1,0],\\
& u(t)=0,\quad t\in[-2,0[,
\end{split}
\end{equation}
which is a particular case of our delayed non-linear optimal control problem (OCP$_\text{D}$)
with $n=m=1$, $a=0$, $b=3$, $r=1$, $s=2$, $g^0\big (x(3)\big )=0$, $f^0(t,x,y,u,v)=x^2+u^2$
and $f(t,x,y,u,v)=yv$. In \cite{Gollmann}, necessary optimality conditions
were proved and applied to \eqref{example1}. The following candidate 
$\big (x^*(\cdot),u^*(\cdot)\big )$ was found:
\begin{equation}
\label{eq:ex:cand:x}
\begin{split}
x^*(t)=
\begin{cases}
1,\ & t\in[-1,2],\\[0.8em]
\displaystyle\frac{e^{t-2}+e^{4-t}}{e^2+1},\ & t\in\ [2,3],
\end{cases}
\end{split}
\end{equation}
and
\begin{equation}
\label{eq:ex:cand:u}
\begin{split}
u^*(t)=
\begin{cases}
0,\ & t\in[-2,0[,\\[0.8em]
\displaystyle\frac{e^t-e^{2-t}}{e^2+1},\ & t\in[0,1],\\[0.8em]
0,\ & t\in\ [1,3].
\end{cases}
\end{split}
\end{equation}
It remains missing in \cite{Gollmann}, however, a proof that such candidate 
\eqref{eq:ex:cand:x}--\eqref{eq:ex:cand:u}
is a solution to the problem. It follows from our
sufficient optimality condition that such claim is indeed true. 
	
We denote that $x_0^*(t)=1$, $t\in[-1,0]$; $x_1^*(t)=1$, 
$t\in[0,1]$; $x_2^*(t)=1$, $t\in[1,2]$; $x_3^*(t)
=\frac{e^{t-2}+e^{4-t}}{e^2+1}$, $t\in[2,3]$; $u_0^*(t)=0$, 
$t\in[-2,0[$; $u_1^*(t)=\frac{e^t-e^{2-t}}{e^2+1}$, $t\in[0,1]$; 
$u_2^*(t)=0$, $t\in[1,2]$ and $u_3^*(t)=0$, $t\in[2,3]$.
Furthermore, the corresponding adjoint function is given by
\begin{equation*}
\begin{split}
\eta(t)=&
\begin{cases}
\eta_1(t),\ & t\in[0,1]\\[0.8em]
\eta_2(t),\ & t\in[1,2]\\[0.8em]
\eta_3(t),\ & t\in[2,3]
\end{cases}\\
=&\begin{cases}
-2t+5+\dfrac{2\big(e^2-1\big)}{\big(e^2+1\big)^2},\ & t\in[0,1]\\[1.4em]
-\left(\dfrac{4e^2}{\big(e^2+1\big)^2}+2\right)t+\dfrac{4\big(e^2
-1\big)}{\big(e^2+1\big)^2}+6+\dfrac{e^{2t-2}-e^{6-2t}}{\big(e^2+1\big)^2},
\ & t\in[1,2]\\[1.4em]
\dfrac{2\big(e^{4-t}-e^{t-2}\big)}{e^2+1},\ & t\in[2,3].
\end{cases}
\end{split}
\end{equation*}
From now on, we are going to ensure that these functions 
satisfy the sufficient optimality conditions studied in 
this section (see Theorem~\ref{Cap4:theo_suf_NLD}). 
So, for $t\in[0,3]$, we intend to find a function $S(t,x)$ 
that is a solution of equation~\eqref{eq_HJ_main_theo} 
with $S\big (3,x(3)\big )=0$. 
As $\eta(t)=\partial_2S\big (t,x(t)\big )$, we obtain that
\begin{equation*}
\begin{split}
S(t,x)=
\begin{cases}
\eta_1(t)x+c_1(t),\ & t\in[0,1]\\[0.8em]
\eta_2(t)x+c_2(t),\ & t\in[1,2]\\[0.8em]
\eta_3(t)x+c_3(t),\ & t\in[2,3],
\end{cases}
\end{split}
\end{equation*}
where $c_i(\cdot)$ is a real function of real variable, $i\in\{1,2,3\}$.
For $t\in[2,3]$, the equation~\eqref{eq_HJ_main_theo} implies that
\begin{align}
\label{ode_c3}
&\ \dot{\eta_3}(t)x^*(t)+\dot{c_3}(t)-\big (x^{*2}(t)
+u^{*2}(t)\big )+\eta_3(t)x^*(t-1)u^*(t-2)=0\nonumber\\
\Leftrightarrow&\ \dot{\eta_3}(t)x_3^*(t)+\dot{c_3}(t)
-\big (x_3^{*2}(t)+u_3^{*2}(t)\big )+\eta_3(t)x_2^*(t-1)u_1^*(t-2)=0\nonumber\\
\Leftrightarrow&-\frac{2\left(e^{4-t}+e^{t-2}\right)}{e^2+1}
\times\frac{e^{t-2}+e^{4-t}}{e^2+1}+\dot{c_3}(t)
-\left(\frac{e^{t-2}+e^{4-t}}{e^2+1}\right)^2\nonumber\\
&+\eta_3(t)\times1\times\frac{e^{t-2}-e^{2-(t-2)}}{e^2+1}=0\nonumber\\
\Leftrightarrow&\ \dot{c_3}(t)=\frac{5\big(e^{2t-4}+e^{8-2t}\big)+2e^2}{\big(e^2+1\big)^2}
\end{align}
with $S\big (3,x(3)\big )=c_3(3)=0$. Solving the differential 
equation~\eqref{ode_c3} with final condition $c_3(3)=0$, 
we obtain that
\begin{equation*}
c_3(t)=\frac{4e^2(t-3)+5\big(e^{2t-4}-e^{8-2t}\big)}{2\big(e^2+1\big)^2}.
\end{equation*}
For $t\in[1,2]$, the equation~\eqref{eq_HJ_main_theo} implies that
\begin{align}
\label{ode_c2}
&\ \dot{\eta_2}(t)x^*(t)+\dot{c_2}(t)-\big (x^{*2}(t)+u^{*2}(t)\big )
+\eta_2(t)x^*(t-1)u^*(t-2)=0\nonumber\\
\Leftrightarrow&\ \dot{\eta_2}(t)x_2^*(t)+\dot{c_2}(t)
-\big (x_2^{*2}(t)+u_2^{*2}(t)\big )+\eta_2(t)x_1^*(t-1)u_0^*(t-2)=0\nonumber\\
\Leftrightarrow&\ -\left(\frac{4e^2}{\big(e^2+1\big)^2}+2\right)
+\frac{2\big(e^{2t-2}+e^{6-2t}\big)}{\big(e^2
+1\big)^2}+\dot{c_2}(t)-1+\eta_2(t)\times1\times0=0\nonumber\\
\Leftrightarrow&\ \dot{c_2}(t)=-\dfrac{2\big(e^{2t-2}+e^{6-2t}
-5e^2\big)-3\big(e^4+1\big)}{\big(e^2+1\big)^2}
\end{align}
with $\eta_2(2)x_2^*(2)+c_2(2)=\eta_3(2)x_3^*(2)+c_3(2)$, 
because $S\big (t,x(t)\big)\in\mathcal{C}^2\left(\mathbb{R}^2,\mathbb{R}\right)$. 
Therefore, the previous condition is equivalent to
\begin{equation}
\label{cond_c2}
c_2(2)=c_3(2)=\frac{5\big(1-e^4\big)-4e^2}{2\big(e^2+1\big)^2}.
\end{equation}
Solving the differential equation~\eqref{ode_c2} with the condition~\eqref{cond_c2}, we have that
\begin{equation*}
c_2(t)=\frac{2t\big(3e^4+10e^2+3\big)+2\big(e^{6-2t}-e^{2t-2}\big)-17e^4-44e^2-7}{2\big(e^2+1\big)^2}.
\end{equation*}
For $t\in[0,1]$, the equation~\eqref{eq_HJ_main_theo} implies that
\begin{align}
\label{ode_c1}
&\ \dot{\eta_1}(t)x^*(t)+\dot{c_1}(t)-\big (x^{*2}(t)+u^{*2}(t)\big)
+\eta_1(t)x^*(t-1)u^*(t-2)=0\nonumber\\
\Leftrightarrow&\ \dot{\eta_1}(t)x_1^*(t)+\dot{c_1}(t)
-\big (x_1^{*2}(t)+u_1^{*2}(t)\big )+\eta_1(t)x_0^*(t-1)u_0^*(t-2)=0\nonumber\\
\Leftrightarrow&-2+\dot{c_1}(t)-1-\left(\frac{e^t-e^{2-t}}{e^2+1}\right)^2
+\eta_1(t)\times1\times0=0\nonumber\\
\Leftrightarrow&\ \dot{c_1}(t)=\frac{e^{4-2t}+e^{2t}+3e^4+4e^2+3}{\big(e^2+1\big)^2}
\end{align}
with $\eta_1(1)x_1^*(1)+c_1(1)=\eta_2(1)x_2^*(1)+c_2(1)$, 
because $S\big (t,x(t)\big )\in\mathcal{C}^2\left(\mathbb{R}^2,\mathbb{R}\right)$. 
Therefore, the previous condition is equivalent to
\begin{equation}
\label{cond_c1}
c_1(1)=c_2(1)=\frac{-9e^4-24e^2-3}{2\big(e^2+1\big)^2}.
\end{equation}
Solving the differential equation~\eqref{ode_c1} 
with the condition~\eqref{cond_c1}, we obtain that
\begin{equation*}
c_1(t)=\frac{2t\big(3e^4+4e^2+3\big)+e^{2t}-e^{4-2t}
-15e^4-32e^2-9}{2\big(e^2+1\big)^2}.
\end{equation*}
	
Concluding, the previous computations show the following result.
	
\begin{prop}
\label{eq:prop1}
Function
\begin{equation*}
\begin{split}
S(t,x)=
\begin{cases}
\eta_1(t)x+c_1(t),\ & t\in[0,1],\\[0.8em]
\eta_2(t)x+c_2(t),\ & t\in[1,2],\\[0.8em]
\eta_3(t)x+c_3(t),\ & t\in[2,3],
\end{cases}
\end{split}
\end{equation*}
with
\begin{equation*}
\begin{split}
\eta_1(t) &= -2t+5+\dfrac{2\big(e^2-1\big)}{\big(e^2+1\big)^2},\\
\eta_2(t) &= -\left(\dfrac{4e^2}{\big(e^2+1\big)^2}
+2\right)t+\dfrac{4\big(e^2-1\big)}{\big(e^2+1\big)^2}
+6+\dfrac{e^{2t-2}-e^{6-2t}}{\big(e^2+1\big)^2},\\
\eta_3(t) &= \dfrac{2\big(e^{4-t}-e^{t-2}\big)}{e^2+1},
\end{split}
\end{equation*}
and
\begin{equation*}
\begin{split}
c_1(t)&=\frac{2t\big(3e^4+4e^2+3\big)+e^{2t}-e^{4-2t}
-15e^4-32e^2-9}{2\big(e^2+1\big)^2},\\
c_2(t)&=\frac{2t\big(3e^4+10e^2+3\big)+2\big(e^{6-2t}
-e^{2t-2}\big)-17e^4-44e^2-7}{2\big(e^2+1\big)^2},\\
c_3(t)&=\frac{4e^2(t-3)+5\big(e^{2t-4}-e^{8-2t}\big)}{2\big(e^2+1\big)^2},
\end{split}
\end{equation*}
is solution of the Hamilton--Jacobi equation~\eqref{eq_HJ_main_theo} 
with $S\big (3,x^*(3)\big )=0$.
\end{prop}


\section{Conclusion}
\label{conclusion}

In this work we did a detailed state of the art associated with optimality conditions 
for delayed optimal control problems. Our survey ends with
sufficient optimality conditions for two different types 
of delayed optimal control problems that are, to the best of our knowledge, 
the first to give an answer to a long-standing open question.
Since the proofs are long, technical, and can be found in \cite{Lemos-Paiao1,Lemos-Paiao2}, 
we did not present them here. However, examples were provided with the purpose 
to illustrate the usefulness of Theorems~\ref{Cap4:theo_OCP-LD} 
and \ref{Cap4:theo_suf_NLD}. As future work, we plan to show the usefulness
of our results to control infectious diseases.


\section*{Acknowledgment}

The authors are strongly grateful 
to the anonymous reviewers for their suggestions 
and invaluable comments.



\end{document}